\documentclass[12pt]{article}
\usepackage{amssymb,color}
\usepackage{algorithmic}
\usepackage{algorithm}

\makeatletter

\@addtoreset{equation}{section}

\makeatother

\newtheorem{thm}[equation]{Theorem}
\newtheorem{lemma}[equation]{Lemma}

\def\P{{\mathbb P}}

\def\Z{{\mathbb Z}}

\def\N{{\mathbb N}}
\def\bX{{\bf X}}
\def\bY{{\bf Y}}

\def\A{{\cal A}}
\def\S{{\cal S}}
\def\sqr{\vcenter{
   \hrule height.1mm
   \hbox{\vrule width.1mm height2.2mm\kern2.18mm\vrule width.1mm}
   \hrule height.1mm}}                  
\def\square{\ifmmode\sqr\else{$\sqr$}\fi}
\def\one{{\bf 1}\hskip-.5mm}
\def\proof{\noindent{\bf Proof. }}

\def\given{\,|\,}

\title{Perfect simulation of a coupling achieving the
 $\bar{d}$-distance between ordered pairs of binary chains of
 infinite order}
\date{September 20, 2010}
\author{Antonio Galves \and Nancy L. Garcia \and Cl\'ementine Prieur}

\begin{document}

\maketitle

\abstract

We explicitly construct a stationary coupling attaining Ornstein's
$\bar{d}$-distance between ordered pairs of binary chains of infinite
order. Our main tool is a representation of the transition
probabilities of the coupled bivariate chain of infinite order as a
countable mixture of Markov transition probabilities of increasing
order. Under suitable conditions on the loss of memory of the chains, this
representation implies that the coupled chain can be represented as a
concatenation of iid sequences of bivariate finite random strings of
symbols. The perfect simulation algorithm is based on the fact that we
can identify the first regeneration point to the left of the origin
almost surely.\\


\section{Introduction}

Let ${\bf X}=(X_n)_{n \in \mathbb{Z}}$ and ${\bf Y}=(Y_n)_{n \in
 \mathbb{Z}}$ be two stationary chains of infinite order on the
alphabet $\A=\{0,1\}$. The $\bar{d}$-distance between ${\bf X}$ and
${\bf Y}$ is defined as
\begin{equation}
\label{eq:dbar}
\overline{d}({\bf X},{\bf Y}) =  \inf \left\{ 
\P(\widetilde{X}_0 \neq \widetilde{Y}_0) \; : \; {(\widetilde{{\bf X}},
\widetilde{{\bf Y}}) \hbox{ stationary coupling of ${\bf X}$ and ${\bf Y}$}} \right\}. 
\end{equation}

The  $\bar{d}$-distance was introduced by Ornstein in several
papers and summarized in an invited article in the first issue of {\it
 The Annals of Probability}  (Ornstein 1973).

The existence of a stationary coupling attaining the $\bar{d}$-minimum follows
from following basic topological considerations. 
\begin{description}
\item (i) The  product space $({\cal A} \times {\cal A})^{\Z}$ is
 compact by Tychonov's Theorem. 
\item (ii) By Prohorov's Theorem, any sequence of probability measures
 on $({\cal A} \times {\cal A})^{\Z}$ has a convergent subsequence in
 the weak$^*$-topology.
\item (iii) Also, the set of all stationary
couplings of ${\bf X}$ and ${\bf Y}$ is a closed subset of the set of
all probability measures on ${\cal A}^{\Z} \times {\cal A}^{\Z}$.
\item (iv) Finally, the Boolean function  $\one_{\{x_0 \neq
 y_0\}}$ that defines the $\bar{d}$-distance is continuous and
bounded.
\end{description}
From (i)--(iv) it follows that there exists at least a coupling which
attains the $\bar{d}$-distance. For more details we refer the reader
to Theorem 4.1 in Villani (2009).     

Obviously this general reasoning does not enable us to explicitly
construct a coupling attaining the $\bar{d}$-minimum.  In spite the
large literature which has been concentrated to this area, as far as
we know the problem of finding explicit solutions was addressed only
for finite alphabet Markov chains and for finite volume Gibbs
measures. To give a further step in this direction is exactly the goal
and the novelty of this paper. We solve in a constructive way the
problem of finding a coupling attaining the $\bar{d}$-distance between
ordered pairs of binary chains of infinite order. First, using basic
stationarity arguments, we prove that the $\bar{d}$-distance is
bounded below by $|\P(Y_0 = 1) - \P(X_0 = 1)|$. Next, we present an
explicit construction of a stationary coupling achieving the infimum
(\ref{eq:dbar}) for stationary chains which are stochastically
ordered. This construction can be effectively implemented in an
algorithmic way to perfectly sample from this minimal
$\bar{d}$-coupling.

This article is organized as follows. In Section \ref{basic} we
introduce the notation and basic definitions. One coupling that
attains the $\bar{d}$-distance is presented in
Section \ref{coupling}.  The perfect sampling algorithm is described
in Section \ref{sec:perfect} and a pseudo-code implementing it is
given by Algorithm \ref{algo1}. The proofs of the
theorems are presented in Sections \ref{proof_decom_regen} and
\ref{proof_main}. We conclude the paper with a final discussion and
some bibliographic remarks (see Section \ref{discuss}).

\section{Basic definitions }\label{basic}

In what follows all the processes and sequences of random variables are
defined on the same probability space $(\Omega, {\cal B}, \P)$. 

Let ${\bf X}=(X_n)_{n \in \mathbb{Z}}$ and ${\bf Y}=(Y_n)_{n \in
 \mathbb{Z}}$ be two stationary chains of infinite order (in the
sense of Harris 1955) on the alphabet $\A=\{0,1\}$. Let $p^X$ and
$p^Y$ respectively be the transition probabilities of these
chains. This means that for any infinite sequence $x^{-1}_{-\infty}
\in \A^{-1}_{-\infty}$ and any symbol $a \in \A$ we have
\[
\P(X_0 = a| X_{-\infty}^{-1} = x_{-\infty}^{-1}) \,=\, p^X(a|x_{-\infty}^{-1})\, ,
\]
\[
\P(Y_0 = a| Y_{-\infty}^{-1} = x_{-\infty}^{-1}) \,=\, p^Y(a|x_{-\infty}^{-1})\, .
\]
In the above formula $x^{-1}_{-\infty}$ denotes the sequence $(x_i)_{i \le
-1}$ and $\A^{-1}_{-\infty}$ the set of all such sequences. These sequences
will be called {\sl pasts}.  Given two integers $m \le n$ we will also use
the notation $x^n_m$ to denote the sequence $(x_m, \ldots,x_n)$ , and
$\A^n_m$ to denote the set of such sequences.

In other terms $p^X$ and $p^Y$ are regular versions of the conditional
expectation of $X_0$ and $Y_0$ with respect to the $\sigma$-algebra
generated by $X^{-1}_{-\infty}$ and $Y^{-1}_{-\infty}$ respectively.

Given two pasts $x^{-1}_{-\infty}$ and $y^{-1}_{-\infty}$, we will say that
$x^{-1}_{-\infty} \le y^{-1}_{-\infty}$, if $x_n \le y_n$ for all $n \le
-1$. This defines a {\sl partial order} on $\A^{-1}_{-\infty}$.

{\bf Condition 1: Ordering condition} We assume that the chains ${\bf
 X}$ and ${\bf Y}$ are stochastically ordered in the following sense
\begin{equation}\label{order}
 p^X(1|x_{-\infty}^{-1})\le  p^Y(1|y_{-\infty}^{-1})\, ,\mbox{\, whenever\,}  x^{-1}_{-\infty} \le y^{-1}_{-\infty}\, .
\end{equation}

The stochastic order between $p^X$ and $p^Y$ makes it possible to construct
a stationary coupling between $\bX$ and $\bY$ in such a way that for all $n
\in \Z$, $X_n \le Y_n$ with probability 1. This coupling is a stationary
chain taking values in the set
\[
\S=\{(0,0), (0,1), (1,1)\}\, .
\]
The transition probabilities $P:\S \times \S_{-\infty}^{-1} \rightarrow [0,
1]$ of this chain are defined as follows: for any pair of ordered pasts
$(x^{-1}_{-\infty} , y^{-1}_{-\infty})\in S_{-\infty}^{-1}$ we have
\begin{eqnarray}
\label{eq:porder}
P\left((1,1)| (x^{-1}_{-\infty} , y^{-1}_{-\infty})\right)&=&
p^X(1|x_{-\infty}^{-1})\, , \nonumber \\
P\left((0,0)| (x^{-1}_{-\infty} , y^{-1}_{-\infty})\right)&=&
p^Y(0|y_{-\infty}^{-1}) \nonumber \\
P\left((0,1)| (x^{-1}_{-\infty} , y^{-1}_{-\infty})\right)&=& p^X(0|x_{-\infty}^{-1})-p^Y(0|y_{-\infty}^{-1})\, .
\end{eqnarray}

We observe that for each pair of ordered pasts $(x^{-1}_{-\infty} ,
y^{-1}_{-\infty})\in S_{-\infty}^{-1}$, $P\left((\cdot,\cdot)|
(x^{-1}_{-\infty} , y^{-1}_{-\infty})\right)$ is the optimal coupling
between $p^X(\cdot |x_{-\infty}^{-1})$ and $p^Y(\cdot |y_{-\infty}^{-1})$.

We want to construct a chain of infinite order on $\S$ invariant with
respect to $P$. This can be done using a regenerative construction of
the chain.  This regenerative construction is based on a decomposition
theorem which states that the stationary chain with infinite memory
can be constructed by choosing at each step, in an iid way, the length
of the suffix of the string of past symbols we need to look in order
to sample the next symbol. 

The above mentioned results will follow under certain conditions on
the transition probabilities:

\noindent {\bf Condition 2: Continuity condition} The transition
probabilities $p^X$ and $p^Y$ on $\A$ are {\it continuous}, that is,
the continuity rates satisfy
$$
\max\left\{\beta^X(k),  \beta^Y(k)\right\}\rightarrow 0 \textrm{ as } k 
\rightarrow \infty\, ,
$$ 
where the continuity rate $\beta^X(k)$ is defined as
\begin{equation}
\label{eq:betak}
\beta^X(k)\,=\,  \max_{a \in \A}
\sup
\{ |p^X(a | x_{-\infty}^{-1})
- p^X(a | y_{-\infty}^{-1})|,
\mbox{ for all } x_{-\infty}^{-1}, \, y_{-\infty}^{-1} \mbox{ with } \,
x_{-k}^{-1} =
y_{-k}^{-1}  \}\, ,
\end{equation}
and similarly for $\beta^Y(k)$. \\

To state our third condition we need some extra notation. For each
pair $(a,b) \in \S$ and each fixed ordered pair of pasts
$(x_{-\infty}^{-1}, y_{-\infty}^{-1}) \in \S_{-\infty}^{-1}$, we
define a non-decreasing sequence
$r_k((a,b)|(x_{-k}^{-1},y_{-k}^{-1}))$ such that

\begin{equation}
r_0(a,b) =  \inf \left\{ P((a,b) \given (u_{-\infty}^{-1},v_{-\infty}^{-1}) ): \, (u_{-\infty}^{-1},v_{-\infty}^{-1}) \in \S_{-\infty}^{-1} \right\} \label{eq:r0} 
\end{equation}
and for $k \ge 1$, $r_k((a,b)|(x_{-k}^{-1},y_{-k}^{-1}))$ is defined by 
\begin{equation}
\inf \left\{ P((a,b) \given (u_{-\infty}^{-1},v_{-\infty}^{-1})):
(u_{-\infty}^{-1},v_{-\infty}^{-1}) \in \S_{-\infty}^{-1},
u_{-k}^{-1}=x_{-k}^{-1}, v_{-k}^{-1}=y_{-k}^{-1}
\right\}\label{eq:rk}
\end{equation}

We then define the non-decreasing sequence $(\alpha_k, k \in \N)$
\begin{equation}
\label{eq:alpha0}
\alpha_0=  \sum_{(a,b) \in \S} r_0((a,b))
\end{equation}
and for $k \ge 1$ 
\begin{equation}
\label{eq:alphak}
\alpha_k((x_{-k}^{-1},y_{-k}^{-1})) = \sum_{(a,b) \in \S} r_k((a,b)|(x_{-k}^{-1},y_{-k}^{-1})),
\end{equation}
and
\begin{equation}
\label{eq:alpha_k}
\alpha_k = \inf \left\{  \alpha_k((x_{-k}^{-1},y_{-k}^{-1})): (x_{-k}^{-1},y_{-k}^{-1}) \in  \S_{-k}^{-1} \right\}.
\end{equation}

\noindent {\bf Condition 3:}  
\begin{equation}\label{alpha}
\prod_{k\ge 0}\alpha_k >0\, .
\end{equation}
\vspace{.3cm}

To better understand Conditions 2 and 3 we will look at an interesting
class of examples which are the renewal processes that forget the past
every time they meet the symbol 1. Take $p^{X}(1|x_{-\infty}^{-1}) =
q^{X}_{\ell(x_{-\infty}^{-1})}$ and $p^{Y}(1|y_{-\infty}^{-1}) =
q^{Y}_{\ell(y_{-\infty}^{-1})}$ where $\ell(u_{-\infty}^{-1}) = \inf
\{ n\ge 1: u_{-n} = 1\}$. We will consider the case when the
expectation of the distance between two successive renewal points in
the $X$ process and the $Y$ process are finite. That is,
$$\sum_{k=0}^{\infty} \prod_{j=0}^{k} (1 - q^{X}_j) < \infty \quad
\mbox{and} \quad \sum_{k=0}^{\infty} \prod_{j=0}^{k} (1 - q^{Y}_j) <
\infty.$$
The divergent case corresponds to the
degenerated case in which the only stationary process with these
transition probabilities is the zero sequence. Notice that
$\sum_{k=0}^{\infty} \prod_{j=0}^{k} (1 - q_j) = \infty $ is weaker
than $\sum_{k=0}^{\infty} q_k < \infty$.

 \paragraph{Example 1:} If $\lim_{k \rightarrow \infty} q^X_k
 \searrow q^X_{\infty} > 0$ and $\lim_{k \rightarrow \infty} q^Y_k
 \searrow q^X_{\infty} >0$ exist, then Condition 2 is satisfied. On
 the other hand, if we take $q_0^X = q^{X}_{2k} \neq q^{X}_{2k+1} =
 q^X_1$ with $0 < q^X_0 < q^X_1 < 1$, Condition 2 is not satisfied.

\paragraph{Example 2:}  If $\lim_{k \rightarrow \infty} q^X_k
 \searrow q^X_{\infty} > 0$ and $\lim_{k \rightarrow \infty} q^Y_k
 \searrow q^X_{\infty} >0$, Condition 3 is equivalent to 
$$ \prod_{n} (1 - q^X_n + q^X_{\infty}) (1 - q^Y_n + q^Y_{\infty}) >
0.$$
For instance, it is enough to have $\sum_n(q_n^X - q_{\infty}^X) =
+\infty$ or  $\sum_n(q_n^Y - q_{\infty}^Y) =
+\infty$ to break Condition 3. 

\section{Construction of our coupling} \label{coupling}

The goal of this section is to present a coupling between the chains
$(X)_n$ and $(Y)_n$ that attains the $\bar{d}$-distance given by
$|\P(Y_0 = 1) - \P(X_0 = 1)|$. To obtain such a coupling Conditions 1--3
are required. Therefore, we assume from now on that they are
satisfied.  

To start the construction we first decompose the transition
probability $P$ given by (\ref{eq:porder}) as a convex combination of
increasing order finite Markov kernels $P_k$ on $\S\times\S_{-k}^{-1}$
for $k \ge 1$. 

Let us define a probability distribution $(\lambda_k, k \in \N)$ as follows.
\begin{equation}
\label{eq:l0}
\lambda_0 = \alpha_0
\end{equation}
and for $k \ge 1$
\begin{equation}
\label{eq:lk}
\lambda_k = \alpha_k \,-\, \alpha_{k-1}.
\end{equation}

The fact that $(\lambda_k, k \in \N)$ is a probability
distribution follows from the fact that $\alpha_k \rightarrow 1$ as
$k$ diverges. Obviously this follows from Condition 2. 

\begin{thm}\label{thm:decomp}  There exists a sequence of transition
 probabilities $P_k$ on $\S\times\S_{-k}^{-1}$ for $k \ge 1$ and a
 probability measure $P_0$ on $\S$ such that for any pair of symbols
 $(a,b)$ in $\S$ and any ordered pair of pasts $(x_{-\infty}^{-1},
 y_{-\infty}^{-1}) \in \S_{-\infty}^{-1}$ we have
\begin{equation}\label{pc}
P\left((a,b)| (x^{-1}_{-\infty} , y^{-1}_{-\infty})\right)
=\lambda_0 P_0\left((a,b)\right) + \sum_{k=1}^{\infty}
\lambda_k P_k\left((a,b)\, |\, (x_{-k}^{-1},y^{-1}_{-k}) \right) \, .
\end{equation}

\end{thm} 

This decomposition allows us to construct simultaneously the pair of
chains $(X_n,Y_n)_{ n \in \Z}$ taking values in $\S$ by concatenating 
bivariate iid strings. This is done as follows.  

Let now ${\bf L}=\{L_n, n \in \Z\}$ be an iid
sequence of random variables such that $\mathbb{P}(L_n=k) = \lambda_k$ where $(\lambda_k, k
\in \N)$ is given by (\ref{eq:l0}) and (\ref{eq:lk}). Define also
$$
T_0 \,=\, \sup \{z \le 0; L_{z+m} \le m, \mbox{ for all } m \ge 0\}
$$
and for $n \ge 1$
$$
T_{-n} \,=\, \sup \{z < T_{-n+1}; L_{z+m} \le m, \mbox{
for all } m \ge 0\}
$$
and
$$
T_{n} \,=\, \inf \{z > T_{n-1}; L_{z+m} \le m, \mbox{
for all } m \ge 0\}.
$$

Given the random variables ${\bf L}=\{L_n, n \in \Z\}$ and ${\bf
 T}=\{T_j, j \in \Z\}$, we construct the bivariate chain
$\{(X_n,Y_n), n \in \Z\}$ by concatenating the bivariate strings
$(X_{T_j}^{T_{j+1}-1},Y_{T_j}^{T_{j+1}-1})$. Each one of these strings
is constructed as follows.

\begin{enumerate}
\item Choose $(X_{T_j},Y_{T_j}) \in \S $ with probability $P_0$ independently of the past.

\item For any $T_j < n \le T_{j+1}-1$ choose  $(X_n,Y_n) \in \S$ with probability $$P_{L_n}\left((\cdot ,\cdot)
| (X^{n-1}_{n-L_n}=x^{n-1}_{n-L_n},Y^{n-1}_{n-L_n}=y^{n-1}_{n-L_n})\right).$$
\end{enumerate}

Observe that if $T_j \le n < T_{j+1}$ then $n - L_n \ge T_j$ and
therefore the choice of the pair $(X_n,Y_n)$ is made independently of
the choice of the symbols
$(X_{-\infty}^{T_j-1},Y_{-\infty}^{T_j-1})$. In this construction, the
transition probabilities $P_k$ are those appearing in Expression
(\ref{pc}).

The existence of infinitely many finite renewal points $T_n$ is given
in the next theorem. 

\begin{thm} \label{thm:regen} The sequence of random times ${\bf T} =
 (T_n, n \in \Z)$ with $\ldots , T_{-1} < T_0 \leq 0 <T_1 < T_2 <
 \ldots$ satisfies \\
 (i) $\mathbb{P}$-almost surely, all the random times $\ldots T_{-1}
 < T_0 \leq 0  <T_1 < T_2  < \ldots$ are finite. \\
 (ii) The random pairs of strings
 $(X_{T_i}^{T_{i+1}-1},Y_{T_i}^{T_{i+1}-1}), \ i \neq 0$ are
 mutually independent and identically distributed. The pair of strings
  $(X_{T_0}^{T_{1}-1},Y_{T_0}^{T_{1}-1})$ is independent from the
  others. 
\end{thm}

We can now present a stationary coupling attaining the $\bar{d}$-distance. This
coupling is obtained concatenating the independent strings
$(X_{T_i}^{T_{i+1}-1},Y_{T_i}^{T_{i+1}-1}), \ i \in \mathbb{Z}$. For this coupling we have the following theorem.

\begin{thm}\label{main} The coupling obtained by concatenating the independent
 strings $(X_{T_i}^{T_{i+1}-1},Y_{T_i}^{T_{i+1}-1}), \ i \in
 \mathbb{Z}$ attains the $\bar{d}$-distance between $\bX$ and $\bY$.

\end{thm}
\section{Perfect simulation algorithm} \label{sec:perfect}

Given two fixed times $m \le n$, we want to perfectly sample
$(X_m^n, Y_m^n)$ according to our minimal $\bar{d}$-coupling between
the chains ${\bf X}$ and ${\bf Y}$ described in Section
\ref{coupling}. 

There is an obvious difficulty: we cannot identify a regeneration
point experimentally. This follows from the fact that, for any
$j \in \Z$ the event ``$j$ is a regeneration point'' is measurable with
respect to the $\sigma$-algebra generated by the random variables
$L_{j+k}, k \ge 0$.

This difficult will be overcome by Algorithm \ref{algo1} whose
pseudo-code is given below. Algorithm \ref{algo1} will produce a
sequence $(\tilde{X}_m^n,\tilde{Y}_m^n)$ as follows. We sequentially
choose iid random variables $L_s, s=n, n-1, \ldots,$ with distribution
$\P(L_s =k) = \lambda_k$ given by (\ref{eq:l0}) and (\ref{eq:lk}).
The algorithm checks every time $t \le m$, until it finds the first
one which has the property that
\[
L_s \le s-t \mbox{\, , for all\, }\, s=t,\ldots, n\, .
\]

Call $T[m,n]$ the first $t \le m$ which has this property:
$$T[m,n] = \sup \{ t \le m; L_s \le s-t, \mbox{ for all } s= t,
\ldots, n\}. 
$$
The random time $T[m,n]$ indicates how far back into the past we have
to look in order to construct $(\tilde{X}_m^n, \tilde{Y}_m^n)$. 

In other terms, if $T[m,n]=t$ then we can choose $(\tilde{X}_t,
\tilde{Y}_t)$ independently of the past with distribution
$P_0$. Moreover, the next pair $(\tilde{X}_{t+1},\tilde{Y}_{t+1})$ can
be chosen using distribution either $P_0$ or
$P_1(\cdot|(\tilde{X}_t,\tilde{Y}_t))$ and recursively we can choose
all the sequence $(\tilde{X}_t^n, \tilde{Y}_t^n)$ without knowledge of
the symbols occurring before time $T[m,n]$. The kernels $P_0$ and
$P_k$ are defined as in Theorem \ref{thm:decomp}.

The sequence $(\tilde{X}_m^n,\tilde{Y}_m^n)$ produced by Algorithm 1
in a finite number of steps depends on the particular choice of the
random variables $L_j, j= T[m,n], \ldots, n$. Let us call this choice
$\tilde{l}_j, j= T[m,n], \ldots, n$. On the other hand, the sequence
$(X_m^n,Y_m^n)$ produced by the theoretical construction presented in
Section \ref{coupling} depends on the choice of $L_j, j \in \Z$. Let
us call $l_j, j \in \Z$ this choice. The important point to stress is
that if $\tilde{l}_j = l_j, j = T[m,n], \ldots, n$ then
$(\tilde{X}_m^n,\tilde{Y}_m^n)=(X_m^n,Y_m^n)$. This is the content of
the following theorem.

We will prove the following theorem.

\begin{thm} \label{thm:perfect} Under Conditions 1--3, for the
  decomposition given by (\ref{pc}), for every pair of integers $m
  \le n$, we have:
\begin{description}
\item (a) $T[m,n]$ is a.s. finite. 
\item (b) The event $\{T[m,n]= t\}$ is measurable with respect to the
 $\sigma$-algebra generated by the random variables $L_s, t \le s \le
 n$. 
\item (c) Algorithm \ref{algo1} stops almost surely after a finite
 number of steps.
\item (d) The sequence $(\tilde{X}_m^n,\tilde{Y}_m^n)$ produced by
 Algorithm 1 is a perfect sample of the minimal $\bar{d}$-coupling
 between the chains ${\bf X}$ and ${\bf Y}$ described in Section
 \ref{coupling}.
\end{description}
\end{thm}

\algsetup{indent=2em}

\begin{algorithm}[!htbp]
\caption{Perfect simulation for a minimal $\bar{d}$-coupling} 
\label{algo1}

\vspace{0.5em}

\begin{algorithmic}[1]
 \REQUIRE Two integers $m \le n$.  \ENSURE The bivariate string
 $(\tilde{X}_m^n, \tilde{Y}_m^n)$ and the past time $T[m,n]$.

\medskip

\STATE $B \leftarrow \emptyset$ \COMMENT{$B$ is the set of time positions $s$ for which the pair $(\tilde{X}_s,\tilde{Y}_s)$ has already been chosen}
\STATE $t \leftarrow m$
\STATE $s \leftarrow t$
\WHILE {$s \le n$} 
  \IF {$s \notin B$}
  \STATE  choose $L_s$  with distribution $\mathbb{P}(L_s=k) = \lambda_k$ independently of everything
          \IF{ $L_s >s-t$}
                      \STATE $t \leftarrow t-1$
                      \STATE $s \leftarrow t$
            \ENDIF  
     \ELSE
      \STATE  choose $(\tilde{X}_s, \tilde{Y}_s)$ with distribution $P_{L_s}((\cdot, \cdot) | (\tilde{X}_{s-L_s}^{s-1}, \tilde{Y}_{s-L_s}^{s-1}))$
      \STATE $B \leftarrow B \cup \{ s \}$
       \STATE $s \leftarrow s+1$
       \ENDIF
\ENDWHILE
\STATE $T[m,n] \leftarrow t$
\RETURN $(\tilde{X}_m^n, \tilde{Y}_m^n)$, $T[m,n]$
\end{algorithmic}
\end{algorithm}

\section{Proofs of Theorems \ref{thm:decomp} and
\ref{thm:regen}}\label{proof_decom_regen}

{\bf Proof of Theorem \ref{thm:decomp}} 

Before starting the proof let us sketch its main ideas. Given an
ordered pair of ``past'' strings
$(x_{-\infty}^{-1},y_{-\infty}^{-1})$, we want to randomly choose a
new random pair of symbols $(a,b) \in \S$ according to $P(
\cdot|(x_{-\infty}^{-1},y_{-\infty}^{-1}))$. This random choice can be
performed as follows.  First make a partition
$\{I((a,b)|(x_{-\infty}^{-1},y_{-\infty}^{-1})), (a,b) \in \S\}$ of
the interval $[0,1]$ where the length of
$I((a,b)|(x_{-\infty}^{-1},y_{-\infty}^{-1})$ is equal to
$P((a,b)|(x_{-\infty}^{-1},y_{-\infty}^{-1}))$. Then, choose a random
element $\xi$ uniformly distributed in $[0,1]$. If $\xi \in
I((a,b)|(x_{-\infty}^{-1},y_{-\infty}^{-1})$, then choose $(a,b)$ as
the new pair of symbols. It turns out that
$I((a,b)|(x_{-\infty}^{-1},y_{-\infty}^{-1})$ can be decomposed as the
following disjoint union
\begin{equation}
\label{idecomp} 
I((a,b)|(x_{-\infty}^{-1},y_{-\infty}^{-1}) \,=\, I_0((a,b)) \cup \cup_{k
    \ge 1} I_k((a,b)|(x_{-k}^{-1},y_{-k}^{-1} )),
\end{equation}
where the length of $I_0((a,b))$ and
$I_k((a,b)|(x_{-k}^{-1},y_{-k}^{-1} )$ are suitably chosen. Loosely
speaking, the length of the interval
$I_k((a,b)|(x_{-k}^{-1},y_{-k}^{-1} )$ is the smallest probability to
choose $(a,b)$ for any pair of ordered pasts having
$(x_{-k}^{-1},y_{-k}^{-1} )$ as ending sequence. 

We can consider a second different partition of $[0,1]$ by using the
increasing sequence $0< \alpha_0 \le \alpha_1 \le \ldots$. The length
of the $k$th element of this partition is precisely
$\lambda_k$. Loosely speaking, if $\xi$ falls on this interval, then
we only need to look at the last $k$ symbols of the past. 

Formally this is done as follows. Let us define a partition of the
interval $[0,1]$ formed by the disjoint intervals
$$I_{0}((0,0)), I_{0}((0,1)),I_{0}((1,1)),$$
and for $k \ge 1$,
$$I_k((0,0)|(x_{-k}^{-1},y_{-k}^{-1})),
I_k((0,1)|(x_{-k}^{-1},y_{-k}^{-1})),
I_k((1,1)|(x_{-k}^{-1},y_{-k}^{-1})), \ldots$$
disposed in the above order in such a way that the left extreme of one
interval coincides with the right extreme of the precedent. These
intervals have length
\begin{equation}
\label{eq:l0x}
|I_{0}((a,b))| = r_0((a,b))
\end{equation}
and for $k \geq 1$, 
\begin{equation}
\label{eq:lkx}
|I_k((a,b) \given (x_{-k}^{-1}, y_{-k}^{-1}))| =
r_k((a,b)|(x_{-k}^{-1}, y_{-k}^{-1})) -
r_{k-1}((a,b)|(x_{-(k-1)}^{-1}, y_{-(k-1)}^{-1})).
\end{equation}

Notice that the continuity of transition probabilities $p^X$ and $p^Y$
implies that 
\begin{equation}
\label{eq:cont}
r_k((a,b)|(x_{-k}^{-1},y_{-k}^{-1})) \rightarrow P((a,b)|(x_{-\infty}^{-1},y_{-\infty}^{-1}))
\end{equation} 
as $k$ diverges.

By construction,
\begin{equation}
\label{eq:p}
P((a,b)|(x_{-\infty}^{-1},y_{-\infty}^{-1})) =  |I_0((a,b))| + \sum_{k
\ge 1} |I_k((a,b)|(x_{-k}^{-1},y_{-k}^{-1} ))|.
\end{equation}

Therefore, we can simulate
$P((a,b)|(x_{-\infty}^{-1},y_{-\infty}^{-1}))$ by using an auxiliary
random variable $\xi$ uniformly distributed on $[0,1]$ as 
\begin{eqnarray}
P((a,b)|(x_{-\infty}^{-1},y_{-\infty}^{-1})) &=& \P\left( \xi \in I_0((a,b)) \cup \cup_{k
    \ge 1} I_k((a,b)|(x_{-k}^{-1},y_{-k}^{-1} )) \right).
\end{eqnarray}
Observe that the right hand side of this equality can be rewritten as
\begin{equation}
\label{eq:decomp}
\sum_{k \geq 0} \P\left( \xi \in [\alpha_{k-1},\alpha_k)\right)\P \left( \xi
  \in  I_0((a,b)) \cup \bigcup_{j
    \ge 1}  I_j((a,b)|(x_{-j}^{-1},y_{-j}^{-1} )) \given \xi \in
  [\alpha_{k-1},\alpha_k)  \right) 
\end{equation}
where $\alpha_{-1}=0$. \\

By construction,
$$[0,\alpha_k) \cap \bigcup_{j > k} I_j((a,b)|(x_{-j}^{-1},y_{-j}^{-1})) =
\emptyset.$$

In other terms, for each $k$, the conditional probabilities on the right
hand side of (\ref{eq:decomp}) depend on the suffix
$(x_{-k}^{-1},y_{-k}^{-1})$ and not on the remaining terms
$(x_{-\infty}^{-(k+1)},y_{-\infty}^{-(k+1)})$. Moreover, 
$$ \sum_{(a,b) \in \S}  \P\left( \xi
\in I_0((a,b)) \cup \bigcup_{j \ge 1} I_j((a,b)|(x_{-j}^{-1},y_{-j}^{-1})) \given
\xi \in [\alpha_{k-1},\alpha_k) \right) = 1.$$

Therefore, we are entitled to define the order $k$ Markov probability
transitions $P_k$ as
\begin{equation}
\label{eq:pxka}
P_k((a,b)|(x_{-k}^{-1},y_{-k}^{-1})) \,=\,  \P\left( \xi
\in I_0((a,b)) \cup \bigcup_{j \ge 1} I_j((a,b)|(x_{-j}^{-1},y_{-j}^{-1})) \given
\xi \in [\alpha_{k-1},\alpha_k) \right).
\end{equation}

Finally we define the probability distribution $(\lambda_k, k \in
\N)$ as follows.
\begin{equation}
\label{eq:lambda0}
\lambda_0 \,=\,  \P( \xi \in [0,\alpha_0)) \,=\, \alpha_0 
\end{equation}
and for $k \ge 1$ 
\begin{equation}
\label{eq:lambdak}
\lambda_k \,=\,  \P( \xi \in [\alpha_{k-1},\alpha_k)) \,=\, \alpha_k -
\alpha_{k-1}.  
\end{equation}

This concludes the proof. \hfill $\square$ \vspace{2mm}

\noindent {\bf Proof of Theorem \ref{thm:regen}} \label{sec:regen}

Define the event $B_n$ as ``n is a regeneration
point''. Formally,  
\begin{equation}
\label{eq:bn}
B_n \,=\, \bigcap_{m \ge 0} \{ L_{n+m} \le m \}. 
\end{equation}

Observe that
\begin{equation}
\label{eq:tn}
\left(\bigcap_{N \ge 1} \bigcup_{n \ge N} B_n \right) \,\cap\,
\left(\bigcap_{N \le 0} \bigcup_{n \le N} B_n \right) \,=\, \bigcap_{k \ge 1} \left\{T_k <
+\infty\right\} \,\cap\, \bigcap_{k \le 0} \left\{T_k > -\infty\right\}.
\end{equation}
Therefore, the existence of infinitely many regeneration times $T_n$
will follow from the following lemma.

\begin{lemma}\label{lemma:bnio}
Assume that $\alpha=\prod_{j=0}^{+\infty}\alpha_j>0$. Then, for any
$N \in \Z$, 
\[
\P \left( \bigcup_{n=N}^{\infty} B_n \right) \,=\, 1\, .
\]
\end{lemma}
\proof
For any $n \in \Z$  define
\[F_n^0 = \{L_{n} >  0\}\]
and $m \ge 1$
\[
F_n^m=\bigcap_{j=0}^{m-1}\{L_{n+j}\le j\}\cap\{L_{n+m} >  m\}\, .
\]
Define
\[
D_1^N=B_N\, ,
\]
and for $k \ge 2$
\[
D_k^N=\bigcup_{n_{1}=N+1}^{+\infty}\ldots\bigcup_{n_{k-1}=n_{k-2}+1}^{+\infty}
\left(F_{N}^{n_1-N-1}\cap\ldots\cap F_{n_{k-2}}^{n_{k-1}-n_{k-2}-1}\cap B_{n_{k-1}}\right)\, .
\]

How to interpret $F_N^m$? Assume $L_N =0$ and therefore, we can choose
$(X_N,Y_N)$ independently of the past symbols
$(X_{-\infty}^{N-1},Y_{-\infty}^{N-1})$. From this point on, we look
at the values of $L_{N+j}$ and we can choose $(X_{N+j},Y_{N+j})$ using
only the knowledge of $(X_N^{N+j-1},Y_N^{N+j-1})$. This sequence
breaks down at $j=m$, since $L_{N+m} > m$ and therefore, the
choice of $(X_{N+m},Y_{N+m})$ depends on the knowledge of symbols
occurring before time $N$. 

Therefore, $D_k^N$ is the event in which the trials,
described above, starting from time $N$ fail exactly $k-1$ times before
finally we find the starting point of a string which is entirely
independent of the past symbols. Therefore, the events $D_k^N$,
$k=1,2,\ldots$ are disjoint and
\[
\bigcup_{n=N}^{+\infty} B_n=\bigcup_{k=1}^{+\infty} D_k^N\, .
\]
Therefore
\[
\P\left(\bigcup_{n=N}^{+\infty} B_n
\right)=\sum_{k=1}^{+\infty}\P(D_k^N)\, .
\]

Since the random lengths $\{L_n, n \in \Z\}$ are
identically distributed, the probabilities computed above do not
depend on the specific choice of $N$. By definition
\[
\P(D_k^N)=\sum_{n_{1}=N+1}^{+\infty}\ldots\sum_{n_{k-1}=n_{k-2}+1}^{+\infty}
\P(F_{N}^{n_1-N-1}\cap\ldots\cap F_{n_{k-2}}^{n_{k-1}-n_{k-2}-1}\cap B_{n_{k-1}})\, .
\]
Using the independence of $F_{N}^{n_1-N-1},
\ldots,F_{n_{k-1}}^{n_k-n_{k-1}-1}$ and $B_{n_k}$ whenever
$N < n_1 <\ldots<n_k$ we can rewrite the right hand side of the last
expression as
\[
\P(D_k^N)=\sum_{n_{1}=N+1}^{+\infty}\ldots\sum_{n_k=n_{k-1}+1}^{+\infty}
\P(F_{N}^{n_1-N-1})\ldots \P(F_{n_{k-1}}^{n_k-n_{k-1}-1}) \P(B_{n_k})\, .
\]

Since  $L_n, n \in \Z$ are iid random variables with $\P(L_0 \le m) = \alpha_m$,
for any $n$, we have 
\begin{eqnarray*}
\P(B_{n}) &=& \P( \cap_{m \ge 0}  \{ L_{n+m} \le m \} ) \\
         & = & \prod_{m \ge 0} \alpha_m  \,=\, \alpha
\end{eqnarray*}
and 
\[
\sum_{l=n+1}^{+\infty}\P(F_n^{l-n-1})= 1-\alpha\, .
\]
Therefore, for any $k \ge 1$ we have
\[
\P(D_k^N)=\alpha (1-\alpha)^{k-1}
\]
and 
\[
\P\left(\bigcup_{n=N}^{+\infty} B_n
\right)=\sum_{k=1}^{+\infty}\alpha (1-\alpha)^{k-1}=1\, .
\]
This
concludes the proof of the lemma. \hfill $\square$

Lemma \ref{lemma:bnio} and the stationarity of the events $B_n$ imply
that 
$$\P\left( \bigcap_{n=-\infty}^0 B_n^c\right) = 0.$$

Observe that for each $n$, if $B_n$ occurs, then
$(X_n^\infty,Y_n^\infty)$ can be chosen independently from from the
past symbols $(X_{-\infty}^{n-1},Y_{-\infty}^{n-1})$. This concludes
the proof of Theorem \ref{thm:regen}. \hfill $\square$

\section{Proof of Theorems \ref{main} and \ref{thm:perfect}}\label{proof_main}

We begin with a lemma giving a lower bound for the $\bar{d}$-distance
between stationary binary chains. For this lemma we are not assuming
that the chains are ordered. 

\begin{lemma} \label{lowerbound}
Let ${\bf X}=(X_n)_{n \in \mathbb{Z}}$ and ${\bf Y}=(Y_n)_{n \in
\mathbb{Z}}$ be any two stationary  chains on
$\{0,1\}$. Then 
$$\bar{d}({\bf X}, {\bf Y}) \ge |\P(Y_0=1) - \P(X_0 = 1)|\, .$$
\end{lemma}

\proof The set of all stationary chains
$(X'_n,Y'_n)_{n \in \Z}$ taking values on $\{0,1\}^2$ such that $\P(X'_n =
1) = \P(X_n=1)$ and $\P(Y'_n =1) = \P(Y_n=1)$ contains the set of all
stationary couplings between $(X_n)_n$ and $(Y_n)_n$. Therefore, 
$$\inf \left\{ 
\P(\widetilde{X}_0 \neq \widetilde{Y}_0) \; : \; {(\widetilde{{\bf X}},
\widetilde{{\bf Y}}) \hbox{ stationary coupling of ${\bf X}$ and ${\bf Y}$}} \right\}$$ is greater than
$$\inf \left\{  \P(\overline{X}_0 \neq \overline{Y}_0) \hbox{ for all  } (\overline{X}_0, \overline{Y}_0) \hbox{ such
   that }\overline{X}_0 \stackrel{\cal D}{=} X_0 \hbox{ and
 }\overline{Y}_0 \stackrel{\cal D}{=} Y_0\right\} \; .$$ It is a
straightforward computation to check that this
last term reaches its minimum with the following optimal coupling between $X_0$
and $Y_0$. For any $a \in \{0,1\}$, take
\begin{eqnarray*}
\P((X'_0,Y'_0)=(a,a)) &=& \min\{\P(X_0=a), \P(Y_0 = a)\}\, ,\\
\P((X'_0,Y'_0)=(a,1-a)) &=& \P(X_0=a) - \P((X'_0,Y'_0)=(a,a))\, .
\end{eqnarray*}
\hfill \square

Now we show that for ordered binary stationary chains $|\P(Y_0=1) -
\P(X_0 = 1)|$ is also an upper bound for $\bar{d}({\bf X}, {\bf Y})$. 

Consider the coupling obtained by concatenating the independent
strings as described in Section \ref{coupling}.  Theorems
\ref{thm:decomp} and \ref{thm:regen} imply that the process
$(X_n,Y_n)_{n \in \Z}$ taking values in $\S$ is stationary. As a
consequence
\begin{itemize} 
\item the chains $(X_n)_{n \in \Z}$ and $(Y_n)_{n \in \Z}$ constructed
simultaneously by the algorithm are also stationary,
\item $(X_0,Y_0)$ is a coupling of the probabilities $\P( X_0 = \cdot)$ and $\P(Y_0 = \cdot)$,
\item moreover by construction $X_0 \le Y_0$. 
\end{itemize}
There exists a unique optimal coupling between $\P( X_0 = \cdot)$ and $\P(Y_0 = \cdot)$,
satisfying the order condition $X_0 \le Y_0$~: 
\[
\P\{(X_0,Y_0)=(0,0) \}= \P(Y_0 = 0)\, ,
\]
\[
\P\{(X_0,Y_0)=(1,1) \}= \P(X_0 = 1)\, ,
\]
\[
\P\{(X_0,Y_0)=(0,1) \}= \P(X_0 = 0)-\P(Y_0 = 0)\, .
\]
With this coupling we have
\begin{equation}\label{probneq}
 \P\{X_0\neq Y_0\}=\, \P(Y_0=1) -  \P(X_0 = 1).
\end{equation}
Equality (\ref{probneq}) together with Lemma \ref{lowerbound}
concludes the proof of Theorem \ref{main}.  \hfill $\square$ \\

To prove Theorem \ref{thm:perfect} let us assume without loss of
generality that $m=0$. 

Assertion (a) follows from the fact that for any $n \geq 0$, $T[0,n]
\ge T_0$ and by Theorem \ref{thm:regen}, $T_0$ is finite almost surely.

The proof of (b) follows from the definition of $T[0,n]$. 

We want to
prove that the number of steps Algorithm \ref{algo1} makes before
stopping is finite. Observe that for each $t$ between $T[0,n]$ and $0$,
the algorithm must do at most $C(|t|+n)$ steps
\begin{itemize}

\item to check if $ L_s\le s-t$ for any $t\le s \le n$

\item and to assign a value to $X_s$ if this is possible.
\end{itemize}  
In the expression $C(|t|+n)$, $C$ is a fixed positive constant which
bounds above the number of operations we need to perform at each
single step.

Therefore the total number of steps Algorithm \ref{algo1} must do
before it stops is bounded above by
\[
C.\sum_{k=0}^{-T[0,n]} (k+n) = C \left[
 (-T[0,n]+1).n+\frac{-T[0,n](-T[0,n]+1)}{2})\right]\, .
\]
This concludes the proof of (c). 

Finally, to prove (d) let us suppose that for $t \le 0$ we have
\begin{equation} \label{rn}
L_t = 0, \, L_{t+1} \le 1, \, \ldots, L_{n} \le n-t.
\end{equation}

Then, the choice of $(X_t^n,Y_t^n)$, according to the theoretical
construction of Section \ref{coupling}, is independent of $L_s, s <
t$. 

By definition, $T[0,n] = \sup \{t \le 0; L_t = 0, \, L_{t+1} \le 1, \,
\ldots, L_{n} \le n-t\}$. By (a) $T[0,n]$ is almost surely finite. By
construction, if $T[0,n]=t$ then 
$$(\tilde{X}_t^n,\tilde{Y}_t^n)=(X_t^n,Y_t^n).$$
\hfill \square   

\section{Final comments and reference remarks}\label{discuss}

The main contribution of this article is to present an explicit
construction of a stationary coupling between ordered binary
chains of infinite order achieving the minimal
$\bar{d}$-distance. Moreover, we show that this explicit construction
is feasible, in the sense that it can be realized by a perfect
simulation algorithm which stops almost surely after a finite number
of steps.

Theorem \ref{main} can be seen as a generalization to the infinite
volume setting of results of Kirillov {\it et al.} (1989) who show
that the classical coupling introduced by Holley (1974) attains
$\bar{d}$-distance for finite volume Gibbs states. Besides Kirillov
{\it et al.} (1989) the only other constructive results on this field
are Ellis (1976, 1978, 1980a, 1980b) which consider the case of Markov
chains on a finite alphabet. Ours seems to be the first constructive
solution for chains of infinite order. Several challenges lay
ahead. For instance the problem of finding a constructive solution for
non-binary chains and/or non-ordered pairs of chains as well as
infinite volume Gibbs measures.

Our results can be presented as a constructive solution for the
Monge-Kantorovich problem with additive cost function on
$C:\A^{\Z}\times \A^{\Z} \rightarrow [0, 1]$
defined as follows. For any pair of sequences $x_{-\infty}^{+\infty}$ and
$y_{-\infty}^{+\infty}$
\[
C(x_{-\infty}^{+\infty},y_{-\infty}^{+\infty})=\sum_{n \in \Z} c_n|x_n -y_n|\, ,
\]
where $(c_n)_{n \in \Z}$ is a sequence of positive real numbers, with
$\sum_{n\in \Z}c_n=1$. This follows straightforward from the following
observation. 
\begin{eqnarray*}
\lefteqn{d_{MK}({\bf X},{\bf Y})} \\
&= & \inf \left\{ 
 \sum_{n \in \Z} c_n \P(\widetilde{X}_n \neq \widetilde{Y}_n) \; : \;
 {(\widetilde{{\bf X}}, \widetilde{{\bf Y}}) \hbox{ stationary
     coupling of ${\bf X}$ and ${\bf Y}$}} \right\} \\
&=& \inf \left\{ 
 \P(\widetilde{X}_0 \neq \widetilde{Y}_0) \sum_{n \in \Z} c_n  \; : \;
 {(\widetilde{{\bf X}}, \widetilde{{\bf Y}}) \hbox{ stationary
     coupling of ${\bf X}$ and ${\bf Y}$}} \right\} \\
&=& \inf \left\{ 
 \P(\widetilde{X}_0 \neq \widetilde{Y}_0)   \; : \;
 {(\widetilde{{\bf X}}, \widetilde{{\bf Y}}) \hbox{ stationary
     coupling of ${\bf X}$ and ${\bf Y}$}} \right\} \\
&=& \overline{d}({\bf X},{\bf Y}).
\end{eqnarray*} 

The Monge-Kantorovich problem has attracted
lots of attention recently.  However, to the best of our knowledge,
ours are the first results in this direction. The literature on MKP is
very extensive. We let the interested reader to find his way starting
with the classical reference Rachev (1984) up to the last Villani
(2009). 

Chains of infinite order seem to have been first studied by Onicescu
and Mihoc (1935a) who called them \emph{chains with complete
connections} (\emph{cha\^{\i}nes \`a liaisons compl\`etes}). 
The name chains of infinite order was coined by Harris (1955).  We
refer the reader to Iosifescu and Grigorescu (1990) for a presentation
of the classical material.  We refer the reader to Fern\'andez,
Ferrari and Galves (2001) for a self contained presentation of chains
of infinite order including the representation of chains of infinite
order as a countable mixture of finite order Markov chains.

Our Theorem \ref{thm:regen} is an application to pairs of chains of
the results in Comets, Fern\'andez and Ferrari (2002). However, our proof
of the result is new and we believe it is simpler than theirs. The
representation of chains of infinite order as a countable mixture of
Markov chains of increasing order appears explicitly in Kalikow
(1990) and implicitly in Ferrari {\it et al.} (2000) and Comets {\it et
 al.}  (2002). Regeneration schemes for chains of infinite order have
been obtained by Berbee (1987) and by Lalley (1986, 2000).

In the literature, the stochastically order between stochastic chains
we considered here is also called {\it domination}. We refer the
reader to the book of Lindvall (1992) for more on the subject.

To assure that Algorithm 1 stops after a finite number of steps we
need weaker conditions than our Conditions 1 and 3. This follows from
the fact that our Algorithm 1 is inspired by the one proposed in
Comets {\it et al.}(2002) in a different context. For details, we
refer the reader to the original article. However, our goal was to
sample from a minimal $\bar{d}$-coupling. It is an open issue if this
can be done under weaker conditions.

\paragraph{Acknowledgments} Thanks to Claire Chauvin, Pierre Collet, D\`avide
Gabrielli and Roberto Imbuzeiro Oliveira for discussions and critical
remarks.  This research was partially supported by PRONEX/ FAPESP's
Project 03/09930-9, FAPESP/CNRS project 2006/50339-0 and CNPq grants
305447/2008-4 (AG) and 301530/2007-6 (NLG).


\begin{thebibliography}{99}
\small

\bibitem{} Berbee, H. (1987) Chains with infinite connections: uniqueness and
Markov representation.  {\it Probab. Theory Related Fields}, {\bf
  76}(2), 243--253.

\bibitem{CFF} Comets, F., Fern\'andez, R.,  Ferrari,  P. A. (2002).
Processes with Long Memory: Regenerative Construction and Perfect
Simulation. {\it Ann. Appl. Probab.} {\bf 12}(3) 921-943.


\bibitem{E1} Ellis, M. H. (1976) The $\overline d$-distance
between two Markov processes cannot always be attained by a Markov
joining. {\it Israel J. Math.} {\bf 24}(3-4), 269--273. 

\bibitem{E2} Ellis, M. H. (1978) Distances between two-state
Markov processes attainable by Markov
joinings. {\it Trans. Amer. Math. Soc.} {\bf 241}, 129--153.


\bibitem{E3} Ellis, M. H. (1980a) Conditions for attaining $\bar
d$ by a Markovian joining. {\it Ann. Probab.} {\bf 8}(3),
431--440. 


\bibitem{E4} Ellis, M. H. (1980b) On Kamae's conjecture concerning the
 $\bar d$-distance between two-state Markov processes. {\it
   Ann. Probab.} {\bf 8}(2), 372--376. 

\bibitem{FFG} Fern\'andez, R., Ferrari, P.A., Galves, A.
(2001). Coupling, renewal and perfect simulation of chains of
infinite order. {\it Notes for a mini-course presented in Vth
  Brazilian School of Probability, 2001}. Can be downloaded from
{\tt http://www.ime.usp.br/\~{}pablo/abstracts/vebp.html}

\bibitem{FMMN} Ferrari, P.A., Maass, A., Martinez, S., Ney,
 P. (2000). Ces\`{a}ro mean distribution of group automata starting
 from measures with summable decay.  {\it Ergodic Theory and
   Dynamical Systems},{\bf 20}(6), 1657--1670.

\bibitem{} Harris, T. E. (1955)  On chains of infinite order.
{\it Pacific J. Math.}, {\bf 5}, 707--724.

\bibitem{} Holley, R. (1974) Remarks on the FKG Inequalities. {\it Commun. math. Phys.}, {\bf 36}, 227--231.
		
\bibitem{} Iosifescu, M., Grigorescu, S., (1990) {\it Dependence with Complete Connections and its
Applications}, Cambridge University Press, Cambridge.  

\bibitem{} Kalikow, S. (1990) Random Markov processes and uniform
martingales., {\it Isr. J. Math.}, {\bf 71}(1), 33-54.

\bibitem{} Kirillov, A.B., Radulescu, D.C., and Styer,
 D.F. (1989). Vasserstein Distances in Two-States Systems. {\it
   Journal of Statistical Physics}, {\bf 56} (5/6), 931--937.
		
\bibitem{} Lalley, S. P. (1986) Regenerative representation for one-dimensional
Gibbs states. {\it  Ann. Probab.}, {\bf  14}(4), 1262--1271.

\bibitem{} Lalley, S. P.  (2000) Regeneration in one-dimensional Gibbs
states and chains with complete connections. {\it Resenhas}, {\bf
  4}(3), 249--281.

\bibitem{} Lindvall, T. (1992) {\it Lectures on the coupling method.}
 Wiley, New York.

\bibitem{} Ornstein, D. (1973) An application of ergodic theory to
 probability theory. {\it Ann. Probab.}, {\bf 1}(1), 43--65. 

\bibitem{} Onicescu, O. and Mihoc, G., (1935) Sur les cha\^{\i}nes
statistiques, {\it C. R. Acad. Sci. Paris}, {\bf 200}, pp. 511--512.

\bibitem{} Rachev, S. T. (1984)  The Monge-Kantorovich problem on mass
transfer and its applications in stochastics.  (Russian)
{\it Teor. Veroyatnost. i Primenen.}, {\bf 29}(4), 625--653. English
translation: Theory Probab. Appl. {\bf 29}(4), 647--676.

\bibitem{} Villani, C. (2009) {\it Optimal transport.
Old and new.}
Grundlehren der Mathematischen Wissenschaften [Fundamental Principles
of Mathematical Sciences], {\bf 338}. Springer-Verlag, Berlin.

				


\end{thebibliography}
\end{document}